\documentclass[12pt]{amsart}
\usepackage{graphicx,amssymb}
\usepackage{graphicx,tikz}
\usepackage{mathrsfs}
\usepackage{amssymb,amsmath,amsthm,color}
\usepackage{graphicx,mcite}
\usepackage{hyperref}
\usepackage{url}
\usepackage{setspace}
\usepackage{float}

\newtheorem{thm}{Theorem}[section]

\theoremstyle{definition}

\theoremstyle{remark}
\newtheorem{rem}[thm]{Remark}
\numberwithin{equation}{section}
\newcommand{\R}{\mathbb R}

\def\TagOnRight

\def\R {\mathbb{R}}
\newcommand{\be}{\begin{equation}}
\newcommand{\ee}{\end{equation}}
\newcommand{\bea}{\begin{eqnarray}}
\newcommand{\eea}{\end{eqnarray}}
\newcommand{\Bea}{\begin{eqnarray*}}
\newcommand{\Eea}{\end{eqnarray*}}
\newcommand{\bt}{\begin{Theorem}}
\newcommand{\et}{\end{Theorem}}

\newcommand{\bpr}{\begin{Proposition}}
\newcommand{\epr}{\end{Proposition}}

\newcommand{\bl}{\begin{Lemma}}
\newcommand{\el}{\end{Lemma}}
\newcommand{\bi}{\begin{itemize}}
\newcommand{\ei}{\end{itemize}}

\begin{document}
\baselineskip16pt

\title[]{$L^{p}$ estimates for multilinear maximal Bochner--Riesz means and square function}
\author[]{Kalachand Shuin}

\address{Department of Mathematics, Indian Institute of Science, Bengaluru-560012, India}
\email{shuin.k.c@gmail.com, kalachands@iisc.ac.in} 

\date{\today}

\date{\today}

\begin{abstract}
In this article we have investigated $L^{p}$ boundedness of the multilinear maximal Bochner--Riesz means and the corresponding square function. We have exploited the ideas given in \cite{JS} to prove our results. 

\end{abstract}
\maketitle

\section{Introduction} \label{section1}
Since the 19th century, convergence of Fourier series has been a central research topic in real variable theory. The study of Bochner--Riesz means is closely related to the investigation of convergence of Fourier series. Additionally, Bochner--Riesz means are among the most significant Fourier multiplier operators in harmonic analysis. For any real number $\alpha\geq0$ and $R>0$, the Bochner--Riesz means are defined by 

\[B_{R}^{\alpha}(f)(x):=\int_{\mathbb{R}^{n}}\Big(1-\frac{|\xi|^{2}}{R^{2}}\Big)^{\alpha}_{+}\hat{f}(\xi)e^{2\pi\iota x\cdot\xi}~d\xi,\]
where $r_{+}=\max\{r,0\}$. It is not hard to check that  $B^{\alpha}_{R}f(x)\lesssim Mf(x)$ for a.e. $x\in\mathbb{R}^{n}$, when  $\alpha>\frac{n-1}{2}$ and $M$ is the Hardy--Littlewood maximal operator.  The index $\alpha=\frac{n-1}{2}$ is recognized as the critical index of Bochner--Riesz means. The primary focus in this area revolves around the   $L^{p}$ boundedness of $B^{\alpha}_{R}$ for minimal value of $\alpha\geq0$.
 %$L^{p}$ boundedness of Bochner--Riesz means correspond to the investigation of 
 %$L^{p}$ convergence of the Fourier series. 
 We define $\alpha(p)=\max\{n|1/p-1/2|-1/2,0\}$ for $1\leq p<\infty$. Then the well-known Bochner--Riesz conjecture posits that for $1\leq p<\infty$, $p\neq2$, the operator $B^{\alpha}_{R}$ maps $L^{p}(\mathbb{R}^{n})$ into $L^{p}(\mathbb{R}^{n})$, if and only if  $\alpha>\alpha(p)$. This conjecture has been settled in dimension $n=1,2$. For $n\geq3$, there are some partial results, see \cite{Carbery, Fefferman1, Fefferman2, Bourgain, BG, Lee0} and the references there in. The point wise convergence of Fourier series corresponds to the study of maximal Bochner--Riesz means defined by 
\[B_{*}^{\alpha}(f)(x):=\sup_{R>0}\Big|\int_{\mathbb{R}^{n}}\Big(1-\frac{|\xi|^{2}}{R^{2}}\Big)^{\alpha}_{+}\hat{f}(\xi)e^{2\pi\iota x\cdot\xi}~d\xi\Big|.\]
The study of $B_{*}^{\alpha}$ is even more challenging when $\alpha<\frac{n-1}{2}$. See \cite{LW, Tao, WG,Jongchon, Christ} for the recent developements of  $L^{p}$ boundedness of the maximal Bochner--Riesz means. In order to study boundedness of  $B_{*}^{\alpha}$, Stein introduced the corresponding square function $G^{\alpha}$ defined by 
\[G^{\alpha}f(x):=\Big(\int^{\infty}_{0}|\frac{\partial}{\partial R}B^{\alpha+1}f(x)|^2 RdR\Big)^{1/2}.\]
See \cite{LRS, Lee1, Sunouchi} and the references there in  for boundedness of the square function $G^{\alpha}$.

A genuine bilinear analogue of Bochner--Riesz means $B^{\alpha}_{R}$ is the bilinear Bochner--Riesz means defined by 
  
  \[\mathcal{B}_{R}^{\alpha}(f_{1},f_{2})(x):=\int_{\mathbb{R}^{2n}}\Big(1-\frac{|\xi_{1}|^{2}+|\xi_{2}|^{2}}{R^{2}}\Big)^{\alpha}_{+}\hat{f_{1}}(\xi_{1})\hat{f_{2}}(\xi_{2})e^{2\pi\iota x\cdot(\xi_{1}+\xi_{2})}~d\xi_{1}d\xi_{2}.\]
  
The study of point wise convergence of two fold product of Fourier series corresponds to the study of maximal bilinear Bochner--Riesz operator $\mathcal{B}^{\alpha}_{*}$, defined by 
\[\mathcal{B}^{\alpha}_{*}(f_{1},f_{2})(x):=\sup_{R>0}|\mathcal{B}^{\alpha}_{R}(f_{1},f_{2})(x)|.\]

The study of bilinear Bochner--Riesz means was first initiated by Bernicot, Grafakos, Song, and Yan in \cite{BGSY}. The  $L^{p}$ boundedness of bilinear Bochner--Riesz means has since been extensively explored by various authors, as noted in \cite{Heping, Grafakos, He, JL}. Recently, Jeong, Lee, and Vargas made significant advancements by introducing a novel approach in dimensions  $n\geq2$. They achieved this by dominating the bilinear operator with the product of two square functions in a pointwise manner, although their method required the conditions 
 $p_1,p_2\geq2$.  More recently, Jotsaroop and Shrivastava \cite{JS} introduced a new decomposition of the bilinear Bochner--Riesz multiplier 
$(1-\frac{|\xi_1|^2+|\xi_2|^2}{R^2})^{\alpha}_{+}$, and in this approach the condition $p_1,p_2\geq2$ are not required.
 The best possible results in this area are presented in \cite{JS}. 
 Also see \cite{CS} and \cite{JSS} for endpoint results  and weighted estimates at the critical index of  bilinear Bochner--Riesz means, respectively.
\subsection{Notations}
Let $n\geq2$. Define $p_{0}(n)=2+\frac{12}{4n-6-k}$ where $n\equiv k$ mod $3$, $k=0,1,2$. Now denote $\mathfrak{p}_{n}=\min\{p_{0}(n),\frac{2(n+2)}{n}\}$. For $n\geq1$ and $1\leq p_{1},p_{2}\leq\infty$ define
\begin{equation*}
\alpha_*(p_1,p_2)=\begin{cases}\alpha(p_1)+\alpha(p_2)&\textup{when}~~ \mathfrak p_n\leq p_1,p_2\leq\infty;\\\\
\alpha(p_1)+\left(\frac{1-2p_2^{-1}}{1-2(\mathfrak p_n)^{-1}}\right)\alpha(\mathfrak p_n)&	\textup{when}~~ \mathfrak p_n\leq p_1\leq\infty~ \text{and}~2\leq p_2<\mathfrak p_n;\\\\
\left(\frac{1-2p_1^{-1}}{1-2(\mathfrak p_n)^{-1}}\right) \alpha(\mathfrak p_n)+\alpha(p_2)&	\textup{when}~~ 2\leq p_1<\mathfrak p_n~ \text{and}~\mathfrak p_n\leq p_2\leq\infty;\\\\
\left(\frac{2-2p_1^{-1}-2p_2^{-1}}{1-2(\mathfrak p_n)^{-1}}\right)\alpha(\mathfrak p_n)&\textup{when}~~2\leq p_1,p_2 <\mathfrak p_n.
\end{cases}
\end{equation*}

For $m\geq2$, $n\geq1$ and Lebesgue measurable functions $F=(f_1,f_2,\dots,f_m)$, we define the multilinear Hardy--Littlewood maximal operator by 
 \[\mathcal{M}_{m}(F)(x):=\sup_{Q\ni x}\prod^{m}_{i=1}\frac{1}{|Q|}\int_{Q}\Big|f_{i}(y)\Big|~dy,\]
 where $Q$ denotes cube in $\mathbb{R}^{n}$ with sides parallel to the coordinate axes. The notation $A\lesssim B$ means that there exists a positive constant $C$ such that $A\leq CB$.
We state the latest results on bilinear maximal Bochner--Riesz means for $n=1$ and $n\geq 2$ separately.  
\begin{thm}\cite{JS} \label{maintheorem:sqr} Let $n\geq 2$ 
	and $(p_1,p_2,p)$ be a triplet such that $p_1,p_2\geq 2$ and $\frac{1}{p}=\frac{1}{p_1}+\frac{1}{p_2},$ then for $\alpha>\alpha_*(p_1,p_2)$ the bilinear maximal Bochner--Riesz operator $\mathcal{B}^{\alpha}_{*}$ maps $L^{p_1}(\R^n)\times L^{p_2}(\R^n)$ into $ L^{p}(\R^n)$. 
\end{thm}
In dimension $n=1$ the bilinear maximal Bochner--Riesz operator $\mathcal B^{\alpha}_{*}$ satisfies the following boundedness properties. 
\begin{thm}\cite{JS}\label{dim1} Let $n=1$ and $1<p_1,p_2<\infty$ be such that $\frac{1}{p_1}+\frac{1}{p_2}=\frac{1}{p}$. The bilinear maximal Bochner--Riesz operator $\mathcal B^{\alpha}_{*}$ maps $L^{p_1}(\R)\times L^{p_2}(\R)$ into $ L^{p}(\R)$ for each of the following cases. 
	\begin{enumerate}
		\item $p_1,p_2\geq 2$ and $\alpha>0$.
		\item $1<p_1<2, p_2\geq 2$ and $\alpha>\frac{1}{p_1}-\frac{1}{2}$.
		\item $1<p_2<2, p_1\geq 2$ and $\alpha>\frac{1}{p_2}-\frac{1}{2}$.
		\item \label{case4} $1<p_1,p_2<2$ and $\alpha>\frac{1}{p}-1$. 
	\end{enumerate}
\end{thm}

In this article we shall consider the multilinear analogue of the Bochner--Riesz means. Let $m\geq3$ and $f_1,f_2,\dots,f_m$ are Schwartz class functions. Then the multilinear Bochner--Riesz means are defined by 
\begin{eqnarray*}
\mathfrak{B}_{R}^{\alpha}(F)(x):=\int_{\mathbb{R}^{mn}}\left(1-\frac{|\vec{\xi}|^{2}}{R^{2}}\right)^{\alpha}_{+}\prod^{m}_{i=1}\hat{f_{i}}(\xi_{i})e^{2\pi\iota x\cdot\sum^{m}_{i=1}{\xi_i}}~d\vec{\xi},
\end{eqnarray*}
where $|\vec{\xi}|^{2}=\sum^{m}_{i=1}|\xi_{i}|^{2}$ and $d\vec{\xi}=\prod^{m}_{i=1}d\xi_i$. 
We shall also consider the multilinear maximal Bochner--Riesz mean $\mathfrak{B}^{\alpha}_{*}$, defined by
\[\mathfrak{B}^{\alpha}_{*}(F)(x):=\sup_{R>0}|	\mathfrak{B}_{R}^{\alpha}(F)(x)|.\] 
It is not hard to check that the critical index of the multilinear Bochner--Riesz means is $\frac{mn}{2}-\frac{1}{2}$. The study of pointwise convergence of $m$-fold product of Fourier series follows from the $L^{p}$ boundedness of  $\mathfrak{B}^{\alpha}_{*}$, i.e. 

\begin{eqnarray}\label{theboundedness}
\Vert \mathfrak{B}^{\alpha}_{*}(F)\Vert_{L^{p}(\mathbb{R}^{n})}\lesssim \prod^{m}_{i=1}\Vert f_{i}\Vert_{L^{p_{i}}(\mathbb{R}^{n})},
\end{eqnarray}
 for $1\leq p_{1},p_{2},\dots,p_{m}\leq\infty$ satisfying $\sum^{m}_{i=1}\frac{1}{p_{i}}=\frac{1}{p}$. 
We also consider the multilinear Bochner--Riesz square function given by 
\begin{eqnarray*}
    \mathcal{G}^{\alpha}(F)(x)&=&\Big(\int^{\infty}_{0}|\frac{\partial}{\partial R}\mathfrak{B}^{\alpha+1}_{R}(F)(x)|^{2} RdR\Big)^{1/2}\\
    &=&\Big(\int^{\infty}_{0}|\mathcal{K}^{\alpha}_{R}\ast F(x)|^2 \frac{dR}{R}\Big)^{1/2},
\end{eqnarray*}
where \[\mathcal{K}^{\alpha}_{R}(\vec{\xi})=2(\alpha+1)\frac{|\vec{\xi}|^2}{R^2}\Big(1-\frac{|\vec{\xi}|^2}{R^2}\Big)^{\alpha}_{+}.\]
Note that for $m=2$,  $L^{p}$ boundedness of the bilinear Bochner--Riesz square function is studied in \cite{CJSS}. 
 
 %\subsection{Notations}
 
\section{Main results}
We shall first state the following results for the trilinear $(m=3)$ maximal Bochner--Riesz mean and trilinear Bochner--Riesz square function for dimension $n\ge1$.  The estimates of the  multilinear $(m\ge4)$ maximal Bochner--Riesz mean will follow from an induction arguments on the $m$-linearity. We shall state our results for dimension $n\ge2$ and $n=1$ separately.
\begin{thm}\label{higherdimension}
	Let  $m=3$, $n\geq2$ and $2\leq p_{1},p_{2},p_{3}\leq\infty$ with $\sum^{3}_{i=1}\frac{1}{p_{i}}=\frac{1}{p}$. Then the inequality \eqref{theboundedness} holds in the following cases
		\begin{enumerate}
		\item For $\mathfrak{p}_{n}\leq p_{1},p_{2},p_{3}\leq\infty$ and $\alpha>\alpha(p_{1})+\alpha(p_{2})+\alpha(p_{3})+\frac{1}{2}$.
		\item For  $i,j,k\in\{1,2,3\} $ and  $i\neq j\neq k$, $2\leq p_{i}<\mathfrak{p}_{n}$, $\mathfrak{p}_{n}\leq p_{j},p_{k}\leq\infty$ and $\alpha>(\frac{1-2p^{-1}_{i}}{1-2\mathfrak{p}^{-1}_{n}})\alpha(\mathfrak{p}_{n})+\alpha(p_{j})+\alpha(p_{k})+\frac{1}{2}$.
		%\item For $2\leq p_{2}<\mathfrak{p}_{n}$, $\mathfrak{p}_{n}\leq p_{1},p_{3}\leq\infty$ and $\alpha>(\frac{1-2p^{-1}_{2}}{1-2\mathfrak{p}^{-1}_{n}})\alpha(\mathfrak{p}_{n})+\alpha(p_{1})+\alpha(p_{3})+\frac{1}{2}$.
		%\item For $2\leq p_{3}<\mathfrak{p}_{n}$, $\mathfrak{p}_{n}\leq p_{1},p_{2}\leq\infty$ and $\alpha>(\frac{1-2p^{-1}_{3}}{1-2\mathfrak{p}^{-1}_{n}})\alpha(\mathfrak{p}_{n})+\alpha(p_{1})+\alpha(p_{2})+\frac{1}{2}$.
		\item For  $i,j,k\in\{1,2,3\} $ and  $i\neq j\neq k$, $2\leq p_{i},p_{j}<\mathfrak{p}_{n}$, $\mathfrak{p}_{n}\leq p_{k}\leq\infty$ and $\alpha>(\frac{2-2p^{-1}_{i}-2p^{-1}_{j}}{1-2\mathfrak{p}^{-1}_{n}})\alpha(\mathfrak{p}_{n})+\alpha(p_{k})+\frac{1}{2}$.
	%\item For $2\leq p_{1},p_{3}<\mathfrak{p}_{n}$, $\mathfrak{p}_{n}\leq p_{2}\leq\infty$ and $\alpha>(\frac{2-2p^{-1}_{1}-2p^{-1}_{3}}{1-2\mathfrak{p}^{-1}_{n}})\alpha(\mathfrak{p}_{n})+\alpha(p_{2})+\frac{1}{2}$.
	%\item For $2\leq p_{2},p_{3}<\mathfrak{p}_{n}$, $\mathfrak{p}_{n}\leq p_{1}\leq\infty$ and $\alpha>(\frac{2-2p^{-1}_{2}-2p^{-1}_{3}}{1-2\mathfrak{p}^{-1}_{n}})\alpha(\mathfrak{p}_{n})+\alpha(p_{1})+\frac{1}{2}$.
	\item For $2\leq p_{1},p_{2},p_{3}<\mathfrak{p}_{n}$  and $\alpha>(\frac{3-2p^{-1}_{1}-2p^{-1}_{2}-2p^{-1}_{3}}{1-2\mathfrak{p}^{-1}_{n}})\alpha(\mathfrak{p}_{n})+\frac{1}{2}$.
	\end{enumerate}
\end{thm}
\begin{thm}\label{dimension1}
	Let $n=1$ and $1<p_{1},p_{2},p_{3}<\infty$ with $\sum^{3}_{i=1}\frac{1}{p_{i}}=\frac{1}{p}$. Then the inequality \eqref{theboundedness} holds in the following cases 
		\begin{enumerate}
		\item For $p_{1},p_{2},p_{3}\geq2$ and $\alpha>1/2$.
		\item For  $i,j,k\in\{1,2,3\} $ and  $i\neq j\neq k$,  $1<p_{i}<2$, $2\leq p_{j},p_{k}<\infty$ and $\alpha>\frac{1}{p_{i}}$. 
		%\item For $1<p_{2}<2$, $2\leq p_{1},p_{3}\leq\infty$ and $\alpha>\frac{1}{p_{2}}$. 
		%\item For $1<p_{3}<2$, $2\leq p_{1},p_{2}<\infty$ and $\alpha>\frac{1}{p_{3}}$.
		\item For   $i,j,k\in\{1,2,3\} $ and  $i\neq j\neq k$, $1<p_{i},p_{j}<2$, $2\leq p_{k}<\infty$ and $\alpha>\frac{1}{p_{i}}+\frac{1}{p_{j}}-\frac{1}{2}$.  
		%\item For $1<p_{1},p_{3}<2$, $2\leq p_{2}\leq\infty$ and $\alpha>\frac{1}{p_{1}}+\frac{1}{p_{3}}-\frac{1}{2}$.  
		%\item For $1<p_{2},p_{3}<2$, $2\leq p_{1}\leq\infty$ and $\alpha>\frac{1}{p_{2}}+\frac{1}{p_{3}}-\frac{1}{2}$.  
		\item  For $1<p_{1},p_{2},p_{3}<2$ and $\alpha>\frac{1}{p}-1$. 
	\end{enumerate}
\end{thm}
\begin{thm}\label{squarefunction3}
   Let $m=3$, $n\geq2$ and $2\leq p_1,p_2,p_3\leq\infty$ with $\frac{1}{p}=\sum^{3}_{i=1}\frac{1}{p_i}$. Then the operator $\mathcal{G}^{\alpha}$ is bounded from $L^{p_1}(\mathbb{R}^{n})\times L^{p_2}(\mathbb{R}^{n})\times L^{p_3}(\mathbb{R}^{n})\to L^{p}(\mathbb{R}^{n})$ if $\alpha$ satisfies the similar conditions as Theorem \ref{higherdimension}. Moreover, for dimension $n=1$, the operator $\mathcal{G}^{\alpha}$ is bounded from $L^{p_1}(\mathbb{R})\times L^{p_2}(\mathbb{R})\times L^{p_3}(\mathbb{R})\to L^{p}(\mathbb{R})$ for $1<p_1,p_2,p_3<\infty$, if  $\alpha$ satisfies the similar conditions as Theorem \ref{dimension1}.
\end{thm}
The results for multilinear $(m\geq4)$ cases are the following.
\begin{thm}\label{mlinear}
    Let $m\geq4$, $n\geq2$ and $2\leq p_i\leq\infty$ with $\sum^{m}_{i=1}1/p_{i}=1/p$. Then \eqref{theboundedness} holds in the following cases 
    \begin{enumerate}
        \item For $\mathfrak{p}_{n}\leq p_{i}\leq\infty$ for all $i=1,2,\dots,m$ and $\alpha>\sum^{m}_{i=1}\alpha(p_{i})+\frac{m-2}{2}$.
		\item For  $i_1,i_2,\dots,i_m\in\{1,2,\dots,m\}$ with $i_j\neq i_k$ if $k\neq j$, and  $2\leq p_{i_1},\dots,p_{i_l}<\mathfrak{p}_n$, $\mathfrak{p}_{n}\leq p_{i_{l+1}},\dots,p_{i_{m}}\leq\infty$, $1\leq l\leq m-1$, and $\alpha>(\frac{l-2\sum^{l}_{k=1}p^{-1}_{i_k}}{1-2\mathfrak{p}^{-1}_{n}})\alpha(\mathfrak{p}_{n})+\sum^{m}_{k=l+1}\alpha(p_{k})+\frac{m-2}{2}$.
  \item For $2\leq p_i<\mathfrak{p}_{n}$, $i=1,2,\dots,m$ and $\alpha>(\frac{m-2\sum^{l}_{k=1}p^{-1}_{i_k}}{1-2\mathfrak{p}^{-1}_{n}})\alpha(\mathfrak{p}_{n})+\frac{m-2}{2}$. 
    \end{enumerate}
\end{thm}
\begin{thm}\label{mlinear1}
    Let $m\geq4$, $n=1$ and $1< p_i<\infty$ with $\sum^{m}_{i=1}1/p_{i}=1/p$. Then \eqref{theboundedness} holds in the following cases 
    \begin{enumerate}
        \item For $p_i\geq2 $, $i=1,2,\dots,m$ and $\alpha>\frac{m-2}{2}$,
        \item For $i_1,i_2,\dots,i_m\in\{1,2,\dots,m\}$ with $i_j\neq i_k$ if $k\neq j$, and  $1< p_{i_1},\dots,p_{i_l}<2$, $2\leq p_{i_{l+1}},\dots,p_{i_{m}}<\infty$, $1\leq l\leq m-1$, and $\alpha>\sum^{l}_{k=1}1/p_{i_k}+\frac{m-l-2}{2}$,
        \item For $1<p_i<2$ for all $i=1,2,\dots,m$ and $\alpha>1/p-1$.
    \end{enumerate}
\end{thm}
The proof of the Theorems \ref{mlinear} and \ref{mlinear1} will follow by an induction argument on the $m(\geq3)$ linearity, and the induction argument will be clear from the proof of $m=3$ case.
\begin{rem}\label{rem}
Let $m\geq4$, $n\geq2$ and $2\leq p_i\leq\infty$ with $1/p=\sum^{m}_{i=1}1/p_i$. Then the multilinear Bochner--Riesz square function $\mathcal{G}^{\alpha}$ maps $\prod^{m}_{i=1}L^{p_i}(\mathbb{R}^{n})\to L^{p}(\mathbb{R}^{n})$ if $\alpha$ satisfies the similar conditions as Theorem \ref{mlinear}. Moreover, for dimension $n=1$, the operator $\mathcal{G}^{\alpha}$ maps $\prod^{m}_{i=1}L^{p_i}(\mathbb{R})\to L^{p}(\mathbb{R})$ if $\alpha$ satisfies the similar conditions as Theorem \ref{mlinear1}.
\end{rem}

\section{Proof of Theorem \ref{higherdimension}: The decomposition of the multiplier}\label{Decomposition}
\begin{proof}
We shall use the idea of \cite{JS} in order to decompose the multiplier $m_{\alpha}(\vec{\xi})=(1-\frac{|\xi_{1}|^{2}+|\xi_{2}|^{2}+|\xi_{3}|^{2}}{R^{2}})^{\alpha}_{+}$. Consider $\psi\in C^{\infty}_{0}[1/2,2]$ and $\psi_{0}\in C^{\infty}_{0}[0,3/4]$ such that \[\sum_{j\geq2}\psi(2^{j}(1-t))+\psi_{0}(t)=1.\] Then we can write 
\begin{eqnarray*}
m_{\alpha}(\vec{\xi})&=&\sum_{j\geq2}\psi(2^{j}(1-\frac{|\xi_{1}|^{2}}{R^{2}}))\varphi_{R}(\xi_{1})^{\alpha}\left(1-\frac{|\eta_{1}|^{2}\varphi^{-1}_{R}(\xi_{1})}{R^{2}}\right)^{\alpha}_{+}+\psi_{0}(\frac{|\xi_{1}|^{2}}{R^{2}})(1-\frac{|\vec{\xi}|^{2}}{R^{2}})^{\alpha}_{+}\\
&=&\sum_{j\geq2}m_{\alpha,j}(\vec{\xi})+m_{\alpha,0}(\vec{\xi}).
\end{eqnarray*}
In the above we have used the notations $\varphi_{R}(\xi_{1})=(1-\frac{|\xi_{1}|^{2}}{R^{2}})_{+}$ and $\eta_{1}=(\xi_{2},\xi_{3})$. In a similar way we shall denote $\eta_{2}=(\xi_{1},\xi_{3})$ and $\eta_{3}=(\xi_{1},\xi_{2})$. 
We shall denote the multiplier $\sum_{j\geq2}m_{\alpha,j}(\vec{\xi})$ by $m_{\alpha,1}(\vec{\xi})$.
Therefore, the operator $\mathfrak{B}^{\alpha}_{R}$ can be decomposed as $$\mathfrak{B}^{\alpha}_{R}(f_{1},f_{2},f_{3})(x)=\mathfrak{B}^{\alpha}_{R,0}(f_{1},f_{2},f_{3})(x)+\mathfrak{B}^{\alpha}_{R,1}(f_{1},f_{2},f_{3})(x),$$
where $\mathfrak{B}^{\alpha}_{R,0}$ corresponds to the multiplier $m_{\alpha,0}$ and $\mathfrak{B}^{\alpha}_{R,1}$ corresponds to the multiplier $m_{\alpha,1}$.

\section{Boundedness of the operator $\mathfrak{B}^{\alpha}_{*,1}$}\label{B1alpha}
Using Stein and Weiss's identity from \cite{SW} [page 278] we get for $j\geq2$, 
\begin{eqnarray*}
m_{\alpha,j}(\vec{\xi})=c_{\delta,\beta}R^{-2\alpha}\psi(2^{j}(1-\frac{|\xi_{1}|^{2}}{R^{2}}))\int^{R_{j}}_{0}(R^{2}\varphi_{R}(\xi_{1})-t^{2})_{+}^{\beta-1}t^{2\delta+1}(1-\frac{|\eta_{1}|^{2}}{t^{2}})^{\delta}_{+}~dt,
\end{eqnarray*}
for $\beta>1/2$ and $\delta>-1/2$ with $\alpha=\beta+\delta$, $c_{\delta,\beta}=2\frac{\Gamma(\delta+\beta+1)}{\Gamma(\delta+1)\Gamma(\beta)}$ and $R_{j}=R\sqrt{2^{-j+1}}$. Therefore, the operator corresponding to $m_{\alpha,j}$ can be written as 
\[\mathfrak{T}^{\alpha}_{j,R}(f_{1},f_{2},f_{3})(x)=c_{\delta,\beta}R^{-2\alpha}\int^{R_{j}}_{0}B^{R,t}_{j,\beta}f_{1}(x)\mathcal{B}^{\delta}_{t}(f_{2},f_{3})(x)t^{2\delta+1}~dt,\]
where 
\begin{eqnarray*}
    &&B^{R,t}_{j,\beta}f_{1}(x):=\int_{\mathbb{R}^n} \psi\Big(2^{j}(1-\frac{|\xi_1|^2}{R^2})\Big)(R^{2}\varphi_{R}(\xi_1)-t^2)^{\beta-1}_{+} e^{2\pi\iota x\cdot\xi_1} \hat{f_{1}}(\xi_1)d\xi_1,\\
    &&\mathcal{B}^{\delta}_{t}(f_2,f_3)(x):=\int_{\mathbb{R}^{2n}} \Big(1-\frac{|\xi_2|^2+|\xi_3|^2}{R^2}\Big)^{\delta}_{+} e^{2\pi\iota x\cdot(\xi_2+\xi_3)} \hat{f_2}(\xi_2)\hat{f_3}(\xi_3)d\xi_2 d\xi_3.
\end{eqnarray*}
Now applying Cauchy--Schwarz inequality we get
\begin{eqnarray*}
\mathfrak{T}^{\alpha}_{j,R}(f_{1},f_{2},f_{3})(x)\leq c_{\delta,\beta}R^{-2\alpha}\left(\int^{R_{j}}_{0}|B^{R,t}_{j,\beta}f_{1}(x)t^{2\delta+1}|^{2}~dt\right)^{1/2}\left(\int^{R_{j}}_{0}|\mathcal{B}^{\delta}_{t}(f_{2},f_{3})(x)|^{2}~dt\right)^{1/2}.
\end{eqnarray*}
Now making a change of variable $t\rightarrow Rt$ in the first component corresponding to $B^{R,t}_{j,\beta}$ we get 
\begin{eqnarray*}
	\mathfrak{T}^{\alpha}_{j,R}(f_{1},f_{2},f_{3})(x)\lesssim 2^{-j/4}\left(\int^{\sqrt{2^{-j+1}}}_{0}|S^{R,t}_{j,\beta}f_{1}(x)t^{2\delta+1}|^{2}~dt\right)^{1/2}\left(\frac{1}{R_{j}}\int^{R_{j}}_{0}|\mathcal{B}^{\delta}_{t}(f_{2},f_{3})(x)|^{2}~dt\right)^{1/2},
\end{eqnarray*}
where 
\[S^{R,t}_{j,\beta}f_{1}(x)=\int_{\mathbb{R}^{n}}\psi(2^{j}(1-\frac{|\xi_{1}|^{2}}{R^{2}}))(1-\frac{|\xi_{1}|^{2}}{R^{2}}-t^{2})^{\beta-1}_{+}\hat{f_{1}}(\xi_{1})e^{2\pi\iota x\cdot \xi_{1}}~d\xi_{1}.\]
Therefore, we get
\begin{eqnarray*}
&&\Vert \mathfrak{T}^{\alpha}_{j,*}(f_{1},f_{2},f_{3})\Vert_{L^{p}(\mathbb{R}^{n})}\\
&&\lesssim 2^{-j/4}\Vert\sup_{R>0}(\int^{\sqrt{2^{-j+1}}}_{0}|S^{R,t}_{j,\beta}f_{1}(x)t^{2\delta+1}|^{2}~dt)^{1/2}\Vert_{L^{p_{1}}(\mathbb{R}^{n})}\\
&&\times\Vert \sup_{R>0}(\frac{1}{R_{j}}\int^{R_{j}}_{0}|\mathcal{B}^{\delta}_{t}(f_{2},f_{3})(x)|^{2}~dt)^{1/2}\Vert_{L^{\tilde{p_{1}}}(\mathbb{R}^{n})}.
\end{eqnarray*}
In the above we have used the following notations. Let  $1\leq p_1,p_2,p_3\leq\infty$ and $1/p=\sum^{3}_{i=1}1/p_i$, then we denote $\frac{1}{\tilde{p}_{i}}=\frac{1}{p}-\frac{1}{p_{i}}$. 
 Note that boundedness of the operator $f_1\rightarrow \sup_{R>0}(\int^{\sqrt{2^{-j+1}}}_{0}|S^{R,t}_{j,\beta}f_{1}(x)t^{2\delta+1}|^{2}~dt)^{1/2}$ follows using Theorem $5.1$ of \cite{JS}.
 \subsection{Boundedness of the operator $(f_{2},f_{3})\rightarrow \sup_{R>0}(\frac{1}{R_{j}}\int^{R_{j}}_{0}|\mathcal{B}^{\delta}_{t}(f_{2},f_{3})(x)|^{2}~dt)^{1/2}$}\label{f2f3}
 Let $d>1$, then we rewrite 
 \begin{eqnarray}
 \mathcal{B}^{\delta}_{t}(f_{2},f_{3})=\sum_{1\leq k\leq d}(\mathcal{B}^{\delta+k-1}_{t}(f_{2},f_{3})-\mathcal{B}^{\delta+k}_{t}(f_{2},f_{3}))+\mathcal{B}^{\delta+d}_{t}(f_{2},f_{3}),
 \end{eqnarray}
 where $d$ is chosen such that $\delta+d>n-1/2$. Now we can write 
 \begin{eqnarray*}
 (\frac{1}{R_{j}}\int^{R_{j}}_{0}|\mathcal{B}^{\delta}_{t}(f_{2},f_{3})(x)|^{2}~dt)^{1/2}&\leq& \sum_{1\leq k\leq d}(\frac{1}{R_{j}}\int^{R_{j}}_{0}|\mathcal{B}^{\delta+k}_{t}(f_{2},f_{3})(x)-\mathcal{B}^{\delta+k-1}_{t}(f_{2},f_{3})(x)|^{2}~dt)^{1/2}\\
 &+& (\frac{1}{R_{j}}\int^{R_{j}}_{0}|\mathcal{B}^{\delta+d}_{t}(f_{2},f_{3})(x)|^{2}~dt)^{1/2}.
 \end{eqnarray*} 
Note that for $\delta+d>n-1/2$, the bilinear Bochner--Riesz means $\mathcal{B}^{\delta+d}_{t}$ is pointwise dominated by the bilinear Hardy--Littlewood maximal operator $\mathcal{M}_2$. Hence the operator $(f_{2},f_{3})\rightarrow \sup_{R>0}(\frac{1}{R_{j}}\int^{R_{j}}_{0}|\mathcal{B}^{\delta+d}_{t}(f_{2},f_{3})(x)|^{2}~dt)^{1/2}$ is bounded from $L^{p_{2}}\times L^{p_{3}}$ to $L^{\tilde{p}_{1}}$ for all $1<p_{2},p_{3}\leq\infty$ and $\frac{1}{p_{2}}+\frac{1}{p_{3}}=\frac{1}{\tilde{p}_{1}}$. On the other hand, the operator $\sup_{R>0}(\frac{1}{R_{j}}\int^{R_{j}}_{0}|\mathcal{B}^{\delta+k}_{t}(f_{2},f_{3})(x)-\mathcal{B}^{\delta+k-1}_{t}(f_{2},f_{3})(x)|^{2}~dt)^{1/2}$ is dominated by the bilinear Bochner--Riesz square function $(\int^{\infty}_{0}|\mathcal{B}^{\delta+k}_{t}(f_{2},f_{3})(x)-\mathcal{B}^{\delta+k-1}_{t}(f_{2},f_{3})(x)|^{2}~\frac{dt}{t})^{1/2}$ with index $\delta+k-1$. Hence, using Theorem $2.1$ of \cite{CJSS} we get the desired estimates. 
\section{Boundedness of the operator $\mathfrak{B}^{\alpha}_{*,0}$}\label{B0alpha}
In order to prove boundedness of the operator $\mathfrak{B}^{\alpha}_{*,0}$, 
we need to further decompose the multiplier $m_{\alpha,0}$ as follows.
\begin{eqnarray*}
m_{\alpha,0}(\vec{\xi})&=&\psi_{0}(\frac{|\xi_{1}|^{2}}{R^{2}})\Big[\sum_{j\geq2}\psi(2^{j}(1-\frac{|\xi_{2}|^{2}}{R^{2}}))\varphi_{R}(\xi_{2})^{\alpha}(1-\frac{|\eta_{2}|^{2}\varphi^{-1}_{R}(\xi_{2})}{R^{2}})^{\alpha}_{+}+\psi_{0}(\frac{|\xi_{2}|^{2}}{R^{2}})(1-\frac{|\vec{\xi}|^{2}}{R^{2}})^{\alpha}_{+}\Big]\\
&:=&m_{\alpha,0,1}(\vec{\xi})+m_{\alpha,0,0}(\vec{\xi}),
\end{eqnarray*}
where  \[m_{\alpha,0,1}(\vec{\xi})=\psi_{0}(\frac{|\xi_{1}|^{2}}{R^{2}})\sum_{j\geq2}\psi(2^{j}(1-\frac{|\xi_{2}|^{2}}{R^{2}}))\varphi_{R}(\xi_{2})^{\alpha}(1-\frac{|\eta_{2}|^{2}\varphi^{-1}_{R}(\xi_{2})}{R^{2}})^{\alpha}_{+}\] and \[m_{\alpha,0,0}(\vec{\xi})=\psi_{0}(\frac{|\xi_{1}|^{2}}{R^{2}})\psi_{0}(\frac{|\xi_{2}|^{2}}{R^{2}})(1-\frac{|\vec{\xi}|^{2}}{R^{2}})^{\alpha}_{+}.\]
Observe that the multiplier $m_{\alpha,0,1}$ is similar to the multiplier $m_{\alpha,1}$ with an extra term $\psi_{0}(\frac{|\xi_{1}|^{2}}{R^{2}})$.
Now we denote the operator corresponding to the multiplier $m_{\alpha,0,1}$ by $\sum_{j\geq2}\mathfrak{T}^{\alpha,1}_{j,R}$. Hence using Stein and Weiss's identity from \cite{SW} [page 278] and a change of variable $t\rightarrow Rt$ we get 
\begin{eqnarray*}
&&\Vert \mathfrak{T}^{\alpha,1}_{j,*}(f_{1},f_{2},f_{3})\Vert_{L^{p}(\mathbb{R}^{n})}\\
&&\lesssim 2^{-j/4}\Vert\sup_{R>0}(\int^{\sqrt{2^{-j+1}}}_{0}|S^{R,t}_{j,\beta}f_{2}(x)t^{2\delta+1}|^{2}~dt)^{1/2}\Vert_{L^{p_{2}}(\mathbb{R}^{n})}\\
&&\times \Vert \sup_{R>0}(\frac{1}{R_{j}}\int^{R_{j}}_{0}|B^{\psi_0}_{R}\mathcal{B}^{\delta}_{t}(f_{1},f_{3})(x)|^{2}~dt)^{1/2}\Vert_{L^{\tilde{p}_{2}}(\mathbb{R}^{n})},
\end{eqnarray*} 
where $\widehat{B^{\psi_0}_{R}\mathcal{B}^{\delta}_{t}(f_{1},f_{3})}(\xi_{1},\xi_{3})=\psi_{0}(\frac{|\xi_{1}|^{2}}{R^{2}})\widehat{\mathcal{B}^{\delta}_{t}(f_{1},f_{3})}(\xi_{1},\xi_{3})$.
Boundedness of the operator $f_{2}\rightarrow \sup_{R>0}(\int^{\sqrt{2^{-j+1}}}_{0}|S^{R,t}_{j,\beta}f_{2}(x)t^{2\delta+1}|^{2}~dt)^{1/2}$ is already discussed in the above. 
\subsection{Boundedness of the operator $(f_{1},f_{3})\rightarrow \sup_{R>0}(\frac{1}{R_{j}}\int^{R_{j}}_{0}|B^{\psi_0}_{R}\mathcal{B}^{\delta}_{t}(f_{1},f_{3})(x)|^{2}~dt)^{1/2}$}\label{f1f3}
Consider $d>1$ such that $\delta+d>n-\frac{1}{2}$. Then we write 
\begin{eqnarray*}
 B^{\psi_0}_{R}\mathcal{B}^{\delta}_{t}(f_{1},f_{3})=\sum_{1\leq k\leq d}(B^{\psi_0}_{R}\mathcal{B}^{\delta+k-1}_{t}(f_{1},f_{3})-B^{\psi_0}_{R}\mathcal{B}^{\delta+k}_{t}(f_{1},f_{3}))+B^{\psi_0}_{R}\mathcal{B}^{\delta+d}_{t}(f_{1},f_{3}).
\end{eqnarray*}
Now applying Cauchy--Schwarz inequality we get
\begin{eqnarray*}
	&&\Big(\frac{1}{R_{j}}\int^{R_{j}}_{0}|B^{\psi_0}_{R}\mathcal{B}^{\delta}_{t}(f_{1},f_{3})(x)|^{2}dt\Big)^{1/2}\\
	&&\leq \sum_{1\leq k\leq d}\Big(\frac{1}{R_{j}}\int^{R_{j}}_{0}|B^{\psi_0}_{R}\mathcal{B}^{\delta+k}_{t}(f_{2},f_{3})(x)-B^{\psi_0}_{R}\mathcal{B}^{\delta+k-1}_{t}(f_{1},f_{3})(x)|^{2}dt\Big)^{1/2}\\
	&&+ \Big(\frac{1}{R_{j}}\int^{R_{j}}_{0}|B^{\psi_0}_{R}\mathcal{B}^{\delta+d}_{t}(f_{1},f_{3})(x)|^{2}~dt\Big)^{1/2}.
\end{eqnarray*}
Since $\delta+d>n-\frac{1}{2}$, therefore the operator $\mathcal{B}^{\delta+d}_{t}$ is dominated by the bilinear Hardy--Littlewood maximal function $\mathcal{M}_2$. Hence we get the desired estimates. On the other hand, boundedness of \[\sum_{1\leq k\leq d}(\frac{1}{R_{j}}\int^{R_{j}}_{0}|B^{\psi_0}_{R}\mathcal{B}^{\delta+k}_{t}(f_{2},f_{3})(x)-B^{\psi_0}_{R}\mathcal{B}^{\delta+k-1}_{t}(f_{1},f_{3})(x)|^{2}dt)^{1/2}\] follows in a similar way as in the Subsection \ref{f2f3}.
\subsection{Boundedness of operator corresponding to the multiplier $m_{\alpha,0,0}(\vec{\xi})$}\label{m00}
Consider two smooth functions $\psi^{1}_{0}\in C^{\infty}_{0}[0,3/32]$ and $\psi^{2}_{0}\in C^{\infty}_{0}[3/64, 3/4]$ such that $\psi_{0}=\psi^{1}_{0}+\psi^{2}_{0}$. Then we decompose the multiplier $m_{\alpha,0,0}$ as follows 
\begin{eqnarray*}
m_{\alpha,0,0}(\vec{\xi})&=&\psi_{0}(\frac{|\xi_{1}|^{2}}{R^{2}})\psi^{1}_{0}(\frac{|\xi_{2}|^{2}}{R^{2}})(1-\frac{|\vec{\xi}|^{2}}{R^{2}})^{\alpha}_{+}+\psi_{0}(\frac{|\xi_{1}|^{2}}{R^{2}})\psi^{2}_{0}(\frac{|\xi_{2}|^{2}}{R^{2}})(1-\frac{|\vec{\xi}|^{2}}{R^{2}})^{\alpha}_{+}\\
&:=&m^{1}_{\alpha,0,0}(\vec{\xi})+m^{2}_{\alpha,0,0}(\vec{\xi}).
\end{eqnarray*}
Now, we shall further decompose the multiplier $m^{1}_{\alpha,0,0}(\vec{\xi})$ as follows 
\begin{eqnarray*}
m^{1}_{\alpha,0,0}(\vec{\xi})&=&\psi_{0}(\frac{|\xi_{1}|^{2}}{R^{2}})\psi^{1}_{0}(\frac{|\xi_{2}|^{2}}{R^{2}})\sum_{j\geq2}\psi(2^{j}(1-\frac{|\xi_{3}|^{2}}{R^{2}}))(1-\frac{|\xi_{3}|^{2}}{R^{2}})^{\alpha}_{+}(1-\frac{|\eta_{3}|^{2}\varphi^{-1}_{R}(\xi_{3})}{R^{2}})^{\alpha}_{+}\\
&+&\psi_{0}(\frac{|\xi_{1}|^{2}}{R^{2}})\psi^{1}_{0}(\frac{|\xi_{2}|^{2}}{R^{2}})\psi_{0}(\frac{|\xi_{3}|^{2}}{R^{2}})(1-\frac{|\vec{\xi}|}{R^{2}})^{\alpha}_{+}.
\end{eqnarray*}
Let  $\mathfrak{T}^{\alpha,1}_{j,R,0}$  denote the Fourier multiplier operator corresponding to the multiplier \[\psi_{0}(\frac{|\xi_{1}|^{2}}{R^{2}})\psi^{1}_{0}(\frac{|\xi_{2}|^{2}}{R^{2}})\psi(2^{j}(1-\frac{|\xi_{3}|^{2}}{R^{2}}))(1-\frac{|\xi_{3}|^{2}}{R^{2}})^{\alpha}_{+}(1-\frac{|\eta_{3}|^{2}\varphi^{-1}_{R}(\xi_{3})}{R^{2}})^{\alpha}_{+}.\] Then using Stein and Weiss's identity and the  change of variable $t\to Rt$ we get 
  \begin{eqnarray*}
  	&&\Vert \mathfrak{T}^{\alpha,1}_{j,*,0}(f_{1},f_{2},f_{3})\Vert_{L^{p}}\\
  	&&\lesssim 2^{-j/4}\Vert\sup_{R>0}(\int^{\sqrt{2^{-j+1}}}_{0}|S^{R,t}_{j,\beta}f_{3}(x)t^{2\delta+1}|^{2}~dt)^{1/2}\Vert_{L^{p_{3}}}\\
  	&&\times\Vert \sup_{R>0}(\frac{1}{R_{j}}\int^{R_{j}}_{0}|B^{\psi_{0}}_{R}B^{\psi^{1}_{0}}_{R}\mathcal{B}^{\delta}_{t}(f_{1},f_{2})(x)|^{2}~dt)^{1/2}\Vert_{L^{\tilde{p}_{3}}},
  \end{eqnarray*} 
  where $\widehat{B^{\psi_{0}}_{R}B^{\psi^{1}_{0}}_{R}\mathcal{B}^{\delta}_{t}}(f_{1},f_{2})(\xi_{1},\xi_{2})=\psi_{0}(\frac{|\xi_{1}|^{2}}{R^{2}})\psi^{1}_{0}(\frac{|\xi_{2}|^{2}}{R^{2}})\widehat{\mathcal{B}^{\delta}_{t}(f_{1},f_{2})}(\xi_{1},\xi_{2})$. The boundedness of this operator will follow in a similar way as in the Subsection \ref{f1f3}.\\
  Now we shall decompose the multiplier $\psi_{0}(\frac{|\xi_{1}|^{2}}{R^{2}})\psi^{1}_{0}(\frac{|\xi_{2}|^{2}}{R^{2}})\psi_{0}(\frac{|\xi_{3}|^{2}}{R^{2}})(1-\frac{|\vec{\xi}|^2}{R^{2}})^{\alpha}_{+}$ as follows
  \begin{eqnarray*}
  \psi_{0}(\frac{|\xi_{1}|^{2}}{R^{2}})\psi^{1}_{0}(\frac{|\xi_{2}|^{2}}{R^{2}})\psi_{0}(\frac{|\xi_{3}|^{2}}{R^{2}})(1-\frac{|\vec{\xi}|^2}{R^{2}})^{\alpha}_{+}&=&\psi_{0}(\frac{|\xi_{1}|^{2}}{R^{2}})\psi^{1}_{0}(\frac{|\xi_{2}|^{2}}{R^{2}})\psi^{1}_{0}(\frac{|\xi_{3}|^{2}}{R^{2}})(1-\frac{|\vec{\xi}|^2}{R^{2}})^{\alpha}_{+}\\
  &+&\psi_{0}(\frac{|\xi_{1}|^{2}}{R^{2}})\psi^{1}_{0}(\frac{|\xi_{2}|^{2}}{R^{2}})\psi^{2}_{0}(\frac{|\xi_{3}|^{2}}{R^{2}})(1-\frac{|\vec{\xi}|^2}{R^{2}})^{\alpha}_{+}.
  \end{eqnarray*}
Observe that the multiplier $\psi_{0}(\frac{|\xi_{1}|^{2}}{R^{2}})\psi^{1}_{0}(\frac{|\xi_{2}|^{2}}{R^{2}})\psi^{1}_{0}(\frac{|\xi_{3}|^{2}}{R^{2}})(1-\frac{|\vec{\xi}|^2}{R^{2}})^{\alpha}_{+}$ has no singularity as its support is away from the boundary of the unit sphere. Hence the corresponding operator will be dominated by the trilinear Hardy--Littlewood maximal function $\mathcal{M}_3$. 
\subsection{Boundedness of the operator corresponding to the multiplier $\psi_{0}(\frac{|\xi_{1}|^{2}}{R^{2}})\psi^{1}_{0}(\frac{|\xi_{2}|^{2}}{R^{2}})\psi^{2}_{0}(\frac{|\xi_{3}|^{2}}{R^{2}})(1-\frac{|\vec{\xi}|^2}{R^{2}})^{\alpha}_{+}$}\label{m00subcase1}
%\begin{proof}
We shall again use Stein and Weiss's identity to rewrite the multiplier as follows 
\begin{eqnarray*}
&&\psi_{0}(\frac{|\xi_{1}|^{2}}{R^{2}})\psi^{1}_{0}(\frac{|\xi_{2}|^{2}}{R^{2}})\psi^{2}_{0}(\frac{|\xi_{3}|^{2}}{R^{2}})(1-\frac{|\vec{\xi}|^2}{R^{2}})^{\alpha}_{+}\\
&&=c_{\delta,\beta}R^{-2\alpha}\psi_{0}(\frac{|\xi_{1}|^{2}}{R^{2}})\psi^{1}_{0}(\frac{|\xi_{2}|^{2}}{R^{2}})\psi^{2}_{0}(\frac{|\xi_{3}|^{2}}{R^{2}})\int^{uR}_{0}(R^{2}\varphi_{R}(\xi_{3})-t^{2})^{\beta-1}_{+}t^{2\delta+1}(1-\frac{|\eta_{3}|^{2}}{t^{2}})^{\delta}_{+}~dt,
\end{eqnarray*}	
where $u=\sqrt{61/64}$. Now we apply Cauchy--Schwarz inequality and a change of variable $t\rightarrow tR$ in the integral corresponding to $f_{3}$. Hence it reduces to prove boundedness of the following operator 
\begin{eqnarray*}
f_{3}\rightarrow \sup_{R>0}\int^{u}_{0}|H^{\beta}_{R,t}f_{3}(x)t^{2\delta+1}|^{2}~dt,
\end{eqnarray*}  
where $H^{\beta}_{R,t}f_{3}(x)=\int_{\mathbb{R}^{n}}\psi^{2}_{0}(\frac{|\xi_{3}|^{2}}{R^{2}})(1-t^{2}-\frac{|\xi_{3}|^{2}}{R^{2}})^{\beta-1}_{+}\hat{f}_{3}(\xi_{3})e^{2\pi\iota x\cdot\xi_{3}}~d\xi_{3}$.
Note that the operator $H^{\beta}_{R,t}f_{3}$ and $S^{R,t}_{j,\beta}f_{3}$ are similar except the term $\psi^{2}_{0}(\frac{|\xi_{3}|^{2}}{R^{2}})$ in place of $\psi(2^{j}(1-\frac{|\xi_{3}|^{2}}{R^{2}}))$ in the multiplier side. But this change does not create any trouble as the function $\psi^{2}_{0}(\frac{|\xi_{3}|^{2}}{R^{2}})$  is compactly supported bounded function with support away from the origin. Hence the analysis can be applied similar to the ones used in proving estimate $(20)$ and Proposition $5.3$ of \cite{JS}. On the other hand boundedness of the operator 
\begin{eqnarray*}
\sup_{R>0}(\frac{1}{R_{j}}\int^{R_{j}}_{0}|B^{\psi_{0}}_{R}B^{\psi^{1}_{0}}_{R}\mathcal{B}^{\delta}_{t}(f_{1},f_{2})(x)|^{2}~dt)^{1/2}
\end{eqnarray*}
is already proved in the Subsection \ref{f1f3}.
%\end{proof}
  \subsection{Boundedness of the operator corresponding to the multiplier $m^{2}_{\alpha,0,0}$}\label{m200}
  We shall further decompose the multiplier as follows 
  \begin{eqnarray*}
  m^{2}_{\alpha,0,0}(\vec{\xi})&=&\psi_{0}(\frac{|\xi_{1}|^{2}}{R^{2}})\psi^{2}_{0}(\frac{|\xi_{2}|^{2}}{R^{2}})\sum_{j\geq2}\psi(2^{j}(1-\frac{|\xi_{3}|^{2}}{R^{2}}))\varphi_{R}(\xi_{3})^{\alpha}(1-\frac{|\eta_{3}|^{2}\varphi^{-1}_{R}(\xi_{3})}{R^{2}})^{\alpha}_{+}\\
  &+&\psi_{0}(\frac{|\xi_{1}|^{2}}{R^{2}})\psi^{2}_{0}(\frac{|\xi_{2}|^{2}}{R^{2}})\psi_{0}(\frac{|\xi_{3}|^{2}}{R^{2}})(1-\frac{|\vec{\xi}|^{2}}{R^{2}})^{\alpha}_{+}.
  \end{eqnarray*}
Observe that the operators corresponding to the multipliers $\psi_{0}(\frac{|\xi_{1}|^{2}}{R^{2}})\psi^{2}_{0}(\frac{|\xi_{2}|^{2}}{R^{2}})\psi(2^{j}(1-\frac{|\xi_{3}|^{2}}{R^{2}}))\varphi_{R}(\xi_{3})^{\alpha}(1-\frac{|\eta_{3}|^{2}\varphi^{-1}_{R}(\xi_{3})}{R^{2}})^{\alpha}_{+}$ are similar as the operators $\mathfrak{T}^{\alpha,1}_{j,*,0}$. Hence the boundedness follows from the boundedness of $\mathfrak{T}^{\alpha,1}_{j,*,0}$.
Note that the multipliers $\psi_{0}(\frac{|\xi_{1}|^{2}}{R^{2}})\psi^{1}_{0}(\frac{|\xi_{2}|^{2}}{R^{2}})\psi^{2}_{0}(\frac{|\xi_{3}|^{2}}{R^{2}})(1-\frac{|\vec{\xi}|^2}{R^{2}})^{\alpha}_{+}$ and $\psi_{0}(\frac{|\xi_{1}|^{2}}{R^{2}})\psi^{2}_{0}(\frac{|\xi_{2}|^{2}}{R^{2}})\psi_{0}(\frac{|\xi_{3}|^{2}}{R^{2}})(1-\frac{|\vec{\xi}|^{2}}{R^{2}})^{\alpha}_{+}$ are similar with the role of $\xi_{3}$ and $\xi_{2}$ interchanged.
\end{proof}
\section{Boundedness of $\mathfrak{B}^{\alpha}_{*}$ in dimension $n=1$}\label{1dimensioncase}
\begin{proof}
We shall use the same decompositions as in the case of dimensions $n\geq2$. Hence using Stein and Weiss's identity and a change of variable followed by the Cauchy--Schwarz inequality we get 
\begin{eqnarray*}
	&&\Vert \mathfrak{T}^{\alpha}_{j,*}(f_{1},f_{2},f_{3})\Vert_{L^{p}(\mathbb{R}^{n})}\\
	&&\lesssim 2^{-j/4}\Vert\sup_{R>0}(\int^{\sqrt{2^{-j+1}}}_{0}|S^{R,t}_{j,\beta}f_{1}(x)t^{2\delta+1}|^{2}~dt)^{1/2}\Vert_{L^{p_{1}}(\mathbb{R}^{n})}\\
 &&\times\Vert \sup_{R>0}(\frac{1}{R_{j}}\int^{R_{j}}_{0}|\mathcal{B}^{\delta}_{t}(f_{2},f_{3})(x)|^{2}~dt)^{1/2}\Vert_{L^{\tilde{p_{1}}}(\mathbb{R}^{n})}.
\end{eqnarray*}
Now invoking Theorem $6.1$ of \cite{JS} and Theorem $2.2$ of \cite{CJSS} we shall get the desired estimates. The rest of the proof can be completed in a similar way as of dimension $n\geq2$. We skip the details to avoid repetition of the arguments. 
\end{proof}
\section{Proof of Theorem \ref{squarefunction3}}
\begin{proof}
We apply the similar approach as Section $3$ of \cite{JS}. 
Consider $\psi\in C^{\infty}_{0}[1/2,2]$ and $\psi_{0}\in C^{\infty}_{0}[0,3/4]$ such that \[\sum_{j\geq2}\psi(2^{j}(1-t))+\psi_{0}(t)=1.\]
This decomposition allows us to rewrite the Fourier transform of the kernel $\mathcal{K}^{\alpha}_{R}$
as 
\begin{eqnarray*}
    \widehat{\mathcal{K}^{\alpha}_{R}}(\vec{\xi})=\sum_{j\geq2}m^{\alpha}_{j,R}(\vec{\xi})+m^{\alpha}_{0,R}(\vec{\xi}),~~\text{where}~~~\vec{\xi}=(\xi_1,\xi_2,\xi_3)\in \mathbb{R}^{3n},
\end{eqnarray*}
and 
\begin{eqnarray*}
    m^{\alpha}_{0,R}(\vec{\xi})&=&\psi_{0}(\frac{|\xi_1|^2}{R^2})\frac{|\vec{\xi}|^2}{R^2}\Big(1-\frac{|\vec{\xi}|^2}{R^2}\Big)^{\alpha}_{+},\\ 
    m^{\alpha}_{j,R}(\vec{\xi})&=&\psi\Big(2^j(1-\frac{|\xi_1|^2}{R^2})\Big)\frac{|\vec{\xi}|^2}{R^2}(1-\frac{|\xi_1|^2}{R^2})^{\alpha}_{+} \Big(1-\frac{|\eta_1|^2}{R^2}(1-\frac{|\xi_1|^2}{R^2})^{-1}\Big)^{\alpha}_{+}, ~\text{for}~j\geq2.
\end{eqnarray*}
We define the corresponding Fourier multiplier operators as follows 
\begin{eqnarray*}
    \frak{g}^{\alpha}_{j,R}(f_1,f_2,f_3)(x)=\int_{\mathbb{R}^{3n}}m^{\alpha}_{j,R}(\xi_1,\xi_2,\xi_3)\prod^{3}_{i=1}\hat{f}_{i}(\xi_i) e^{2\pi\iota x\cdot(\xi_1+\xi_2+\xi_3)} d\vec{\xi},~\text{for}~j=0,2,3,\dots.
\end{eqnarray*}
Therefore, 
\begin{eqnarray*}
\mathcal{G}^{\alpha}(f_1,f_2,f_3)(x)&\leq& \mathcal{G}_{0}^{\alpha}(f_1,f_2,f_3)(x)+\sum_{j\geq2}\mathcal{G}_{j}^{\alpha}(f_1,f_2,f_3)(x),~\text{where}\\
    \mathcal{G}^{\alpha}_{j}(f_1,f_2,f_3)(x)&=&\Big(\int^{\infty}_{0}|\frak{g}^{\alpha}_{j,R}(f_1,f_2,f_3)(x)|^{2}\frac{dR}{R}\Big)^{1/2}, ~\text{for}~j=0,2,3,\dots 
\end{eqnarray*}
Now it remains to show that
\begin{eqnarray*}
    \Vert \mathcal{G}^{\alpha}_{j}\Vert_{L^{p_1}\times L^{p_2}\times L^{p_3}\to L^{p}}\leq C_j,~\text{s.t.}~ \sum_{j\geq2}C_j<\infty.
\end{eqnarray*}
\subsection{Boundedness of $\mathcal{G}^{\alpha}_{j}$ for $j\geq2$}
We further decompose the operator $\mathfrak{g}^{\alpha}_{j,R}$ using the identity of Stein and Weiss [page 278, \cite{SW}]. We define $\varphi_{R}(\xi_i)=(1-\frac{|\xi_i|^2}{R^2})_{+}$. Therefore, 
\[(1-\frac{|\eta_1|^2}{R^2\varphi_{R}(\xi_1)})^{\alpha}_{+}=c_{\delta,\beta}R^{-2\alpha}\varphi_{R}(\xi_1)^{-\alpha}\int^{R}_{0}(R^2\varphi_{R}(\xi_1)-t^2)^{\beta-1}_{+}t^{2\delta+1}(1-\frac{|\eta_1|^2}{R^2})^{\delta}_{+}~dt,\]
where $\alpha=\beta+\delta$ with $\beta>1/2$, $\delta>-1/2$ and $c_{\delta,\beta}=2\frac{\Gamma(\delta+\beta+1)}{\Gamma(\delta+1)\Gamma(\beta)}$.
Hence we get,
\begin{eqnarray*}
    \mathfrak{g}^{\alpha}_{j,R}(f_1,f_2,f_3)(x)=c_{\delta,\beta}R^{-2\alpha}\int^{R_j}_{0}\Big(B^{R,t}_{j,\beta}f_1(x)\mathcal{B}^{\delta}_{t}(f_2,f_3)(x)+A^{R,t}_{j,\beta}f_1(x)\mathcal{A}^{\delta}_{t}(f_2,f_3)(x)\Big) t^{2\delta+1} dt,
\end{eqnarray*}
where $R_j=R\sqrt{2^{-j+1}}$ and 
\begin{eqnarray*}
    B^{R,t}_{j,\beta}f_1(x)&=&\int_{\mathbb{R}^n}\hat{f_1}(\xi_1)\psi\Big(2^j(1-\frac{|\xi_1|^2}{R^2})\Big)\frac{|\xi_1|^2}{R^2}(R^2\varphi_{R}(\xi_1)-t^2)^{\beta-1}e^{2\pi\iota x\cdot\xi_1} d\xi_1,\\
    A^{R,t}_{j,\beta}f_1(x)&=&\int_{\mathbb{R}^n}\hat{f_1}(\xi_1)\psi\Big(2^j(1-\frac{|\xi_1|^2}{R^2})\Big)(R^2\varphi_{R}(\xi_1)-t^2)^{\beta-1}_{+}e^{2\pi\iota x\cdot\xi_1} d\xi_1,\\
    \mathcal{A}^{\delta}_{t}(f_2,f_3)(x)&=&\int_{\mathbb{R}^{2n}}\hat{f_2}(\xi_2)\hat{f_3}(\xi_3)\frac{|\eta_1|^2}{R^2}(1-\frac{|\eta_1|^2}{t^2})^{\delta}_{+} e^{2\pi\iota x\cdot(\xi_2+\xi_3)} d\xi_2 d\xi_3.
\end{eqnarray*}
Therefore, \[\mathfrak{g}^{\alpha}_{j,R}(f_1,f_2,f_3)(x)=T^1_{j,R}+T^2_{j,R},\] where 
\begin{eqnarray*}
    T^1_{j,R}&=&c_{\delta,\beta}R^{-2\alpha}\int^{R_j}_{0}B^{R,t}_{j,\beta}f_1(x)\mathcal{B}^{\delta}_{t}(f_2,f_3)(x) t^{2\delta+1} dt,\\
    T^2_{j,R}&=&c_{\delta,\beta}R^{-2\alpha}\int^{R_j}_{0}A^{R,t}_{j,\beta}f_1(x)\mathcal{A}^{\delta}_{t}(f_2,f_3)(x) t^{2\delta+1} dt.
\end{eqnarray*}
Now using Cauchy--Schwarz inequality and the ideas of [page $19,20$ of  \cite{CJSS}] we get 
\begin{eqnarray*}
    \Vert \Big(\int^{\infty}_{0}|T^2_{j,R}|^2\frac{dR}{R}\Big)^{1/2}\Vert_{L^p}\lesssim 2^{j(\alpha(p_1)-\beta-1+\epsilon)}\prod^{3}_{i=1}\Vert f_i\Vert_{L^{p_i}},
\end{eqnarray*}
for $p_1\geq\frak{p}_{n}$ or $p_1=2$ when $\beta>\alpha(p_1)+1/2$, and $p_2,p_3\geq \frak{p}_n$ or either of $p_2,p_3=2$ or both $p_2,p_3=2$ when $\delta>\alpha_{*}(p_2,p_3)$, and the exponent of $2^j$ is negative. 
On the other hand, using Cauchy--Schwarz inequality we get 
\begin{eqnarray*}
     &&\Vert \Big(\int^{\infty}_{0}|T^1_{j,R}|^2\frac{dR}{R}\Big)^{1/2}\Vert_{L^p}\\
     &&\lesssim \Vert \Big[\int^{\infty}_{0}\Big(\int^{R_j}_{0}|B^{R,t}_{j,\beta}f_1(x)t^{2\delta+1}|^2 dt\Big)\Big(\int^{R_j}_{0}|\mathcal{B}^{\delta}_{t}(f_2,f_3)(x)|^2 dt\Big) \frac{dR}{R^{4\alpha+1}}\Big]^{1/2}\Vert_{L^p}\\
     &&\lesssim 2^{-j/4} \Vert \Big[\int^{\infty}_{0}\Big(\int^{R_j}_{0}|S^{R,t}_{j,\beta}f_1(x)t^{2\delta+1}|^2 dt\Big) \frac{dR}{R}\Big]^{1/2}\Vert_{L^{p_1}} \Vert \sup_{R>0}\Big(\frac{1}{R_j}\int^{R_j}_{0}|\mathcal{B}^{\delta}_{t}(f_2,f_3)(x)|^2 dt\Big)^{1/2}\Vert_{L^{\tilde{p}_1}},
\end{eqnarray*}
where $1/\tilde{p}_1=1/p_2+1/p_3$ and 
\[S^{R,t}_{j,\beta}f_1(x)=\int \psi\Big(2^j(1-\frac{|\xi_1|^2}{R^2})\Big) \frac{|\xi_1|^2}{R^2} (1-\frac{|\xi_1|^2}{R^2}-t^2)^{\beta-1}_{+} \hat{f_1}(\xi_1)e^{2\pi\iota x\cdot \xi_1} d\xi_1.\]
Using the similar idea as Lemma $4.2$ of \cite{JS} we get that the operator 
\[(f_2,f_3)\to \sup_{R>0}(\frac{1}{R}\int^{R}_{0}|\mathcal{B}^{\delta}_{t}(f_2,f_3)(x)|^2 dt)^{1/2}\]
is bounded from $L^{p_2}\times L^{p_3}\to L^{\tilde{p_1}}$ when $\delta>\alpha_{*}(p_2,p_3)$. Moreover, using Theorem $4.2$ of \cite{CJSS} we get 
\[\Vert (\int^{\infty}_{0}|S^{R,t}_{j,\beta}f_1(x)|^2\frac{dR}{R})^{1/2}\Vert_{L^{p_1}}\lesssim 2^{j(\alpha(p_1)-\beta-\frac{1}{4}+\epsilon)}\Vert f_1\Vert_{L^{p_1}},\]
when $\beta>\alpha(p_1)+1/2$. This yields the desired estimates.
\subsection{Boundedness of $\mathcal{G}^{\alpha}_{0}$}
Boundedness of $\mathcal{G}^{\alpha}_{0}$ can be completed using the similar ideas as the previous section with minor modifications. In fact, the boundedness of $\mathfrak{g}^{\alpha}_{0}$ follows from another decomposition w.r.t the variable $\xi_2$ or $\xi_3$. See Section \ref{B0alpha} for the decomposition part. 
\subsection{Boundedness of $\mathcal{G}^{\alpha}$ in dimension $n=1$}
We apply the same decomposition as $n\geq2$ case and walk the same path till we reduce the operator  $\mathfrak{g}^{\alpha}_{j,R}$ into 
\[\mathfrak{g}^{\alpha}_{j,R}(f_1,f_2,f_3)(x)=T^1_{j,R}+T^{2}_{j,R}.\]
Now in order to handle $T^{2}_{j,R}$ we invoke Theorem $6.1$ of \cite{JS}  to get 
\begin{eqnarray*}
    \Vert \sup_{R>0}\Big(\int^{\sqrt{2^{-j+1}}}_{0}|\tilde{S}^{R,t}_{j,\beta}f_1(\cdot)t^{2\delta+1}|^2 dt\Big)^{1/2}\Vert_{L^{p_1}}\lesssim 2^{j(1/4-\tilde{\alpha})}\Vert f_1\Vert_{L^{p_1}},
\end{eqnarray*}
for $\beta>1/2$, when $2\leq p_1\leq\infty$, and for $\beta>1/p_1$ when $1<p_1<2$. Here $\tilde{\alpha}=\min\{\alpha,\delta+1/2\}$. After a change of variable, the operator $\mathcal{A}^{\delta}_{t}$ satisfies 
\begin{eqnarray*}
    \Vert \Big[\int^{\infty}_{0}\Big(\int^{\sqrt{2^{-j+1}}}_{0}|\mathcal{A}^{\delta}_{Rt}(f_2,f_3)(x)|^2 dt\Big)\frac{dR}{R}\Big]^{1/2}\Vert_{L^{\tilde{p}_1}}\lesssim 2^{-5j/4}\Vert f_2\Vert_{L^{p_2}}\Vert f_3\Vert_{L^{p_3}},
\end{eqnarray*}
for $(p_2,p_3,\delta)$ satisfying the conditions  of $(p_1,p_2,\alpha)$ of Theorem \ref{dim1}. On the other hand,
\begin{eqnarray*}
     &&\Vert (\int^{\infty}_{0}|T^1_{j,R}|^2\frac{dR}{R})^{1/2}\Vert_{L^p}\\
     &&\lesssim 2^{-j/4} \Vert \Big[\int^{\infty}_{0}\Big(\int^{R_j}_{0}|S^{R,t}_{j,\beta}f_1(x)t^{2\delta+1}|^2 dt\Big) \frac{dR}{R}\Big]^{1/2}\Vert_{L^{p_1}} \Vert \sup_{R>0}\Big(\frac{1}{R_j}\int^{R_j}_{0}|\mathcal{B}^{\delta}_{t}(f_2,f_3)(x)|^2 dt\Big)^{1/2}\Vert_{L^{\tilde{p}_1}}.
\end{eqnarray*}
Now from equation $(28)$ of \cite{CJSS} we have 
\[\Vert \Big[\int^{\infty}_{0}\Big(\int^{R_j}_{0}|S^{R,t}_{j,\beta}f_1(x)t^{2\delta+1}|^2 dt\Big) \frac{dR}{R}\Big]^{1/2}\Vert_{L^{p_1}} \lesssim 2^{-(j-1)(\delta+\frac{1}{2})}\Vert f_1\Vert_{L^{p_1}},\]
for $\beta>1/2$ when $2\leq p_1<\infty$ and for $\beta>1/p_1$ when $1<p_1<2$. And for the second term we have 
\[\Vert \sup_{R>0}\Big(\frac{1}{R_j}\int^{R_j}_{0}|\mathcal{B}^{\delta}_{t}(f_2,f_3)(x)|^2 dt\Big)^{1/2}\Vert_{L^{\tilde{p}_1}}\lesssim \Vert f_2\Vert_{L^{p_2}}\Vert f_3\Vert_{L^{p_3}},\]
where $(p_2,p_3,\delta)$ satisfies the conditions of $(p_1,p_2,\alpha)$ of Theorem \ref{dim1}. This completes the proof.
\end{proof}
\section{Proof of Theorems \ref{mlinear} and \ref{mlinear1}}
\begin{proof}
Boundedness of the $m$-linear maximal Bochner--Riesz operator can be proved using the induction argument on the $m$-linearity of the operator. We define the $N$-linear operator $\mathfrak{B}^{\alpha}_{*}$ and  we have already shown boundedness of $\mathfrak{B}^{\alpha}_{*}$ for $N=3$. Also we have proved boundedness of the Bochner--Riesz square function $\mathcal{G}^{\alpha}$ for $N=3$. By the induction hypothesis we assume that Theorem \ref{mlinear} holds for $N=m-1$, and the Remark \ref{rem} holds for $N=m-1$. Then we claim to prove Theorem \ref{mlinear}  and the Remark \ref{rem} for $N=m$.
First we define a few more notations here. Given $\vec{\xi}=(\xi_1,\xi_2,\dots,\xi_m)$ we consider  $\eta_1=(\xi_2,\dots,\xi_m) $, $\eta_m=(\xi_1,\xi_2,\dots,\xi_{m-1})$ and $\eta_j=(\xi_1,\xi_2,\dots,{\xi_{j-1}},\xi_{j+1},\dots, \xi_m)$ for all $j=2,\dots,m-1$. Also define  $F_{1}=(f_2,\dots,f_m)$, $F_m=(f_1,f_2,\dots,f_{m-1})$ and  $F_{j}=(f_1,f_2,\dots,f_{j-1},f_{j+1},\dots,f_m)$ for $j=2,3,\dots,m-1$.
Consider $\psi\in C^{\infty}_{0}[1/2,2]$ and $\psi_{0}\in C^{\infty}_{0}[0,3/4]$ such that \[\sum_{j\geq2}\psi(2^{j}(1-t))+\psi_{0}(t)=1.\] Then we write 
\begin{eqnarray*}
m_{\alpha}(\vec{\xi})&=&\sum_{j\geq2}\psi(2^{j}(1-\frac{|\xi_{1}|^{2}}{R^{2}}))\varphi_{R}(\xi_{1})^{\alpha}\left(1-\frac{|\eta_{1}|^{2}\varphi^{-1}_{R}(\xi_{1})}{R^{2}}\right)^{\alpha}_{+}+\psi_{0}(\frac{|\xi_{1}|^{2}}{R^{2}})(1-\frac{|\vec{\xi}|^{2}}{R^{2}})^{\alpha}_{+}\\
&=&\sum_{j\geq2}m_{\alpha,j}(\vec{\xi})+m_{\alpha,0}(\vec{\xi}).
\end{eqnarray*}
We shall denote the multiplier $\sum_{j\geq2}m_{\alpha,j}(\vec{\xi})$ by $m_{\alpha,1}(\vec{\xi})$.
Therefore, the operator $\mathfrak{B}^{\alpha}_{R}$ can be written as $$\mathfrak{B}^{\alpha}_{R}(F)(x)=\mathfrak{B}^{\alpha}_{R,0}(F)(x)+\mathfrak{B}^{\alpha}_{R,1}(F)(x),$$
where $\mathfrak{B}^{\alpha}_{R,0}$ corresponds to the multiplier $m_{\alpha,0}$ and $\mathfrak{B}^{\alpha}_{R,1}$ corresponds to the multiplier $m_{\alpha,1}$.
\subsection{Boundedness of $\mathfrak{B}^{\alpha}_{*,1}$}
Invoking the identity of Stein and Weiss as in Section \ref{B1alpha} we get 
\[\mathfrak{I}^{\alpha}_{j,R}(F)(x)=c_{\delta,\beta}R^{-2\alpha}\int^{R_{j}}_{0}B^{R,t}_{j,\beta}f_1(x)\mathfrak{B}^{\delta}_{t}(F_1)(x) t^{2\delta+1} dt.\]
Now applying Cauchy--Schwarz inequality  and a change of variable we get 
\begin{eqnarray*}
\mathfrak{T}^{\alpha}_{j,R}(F)(x)\lesssim 2^{-j/4}\left(\int^{\sqrt{2^{-j+1}}}_{0}|S^{R,t}_{j,\beta}f_{1}(x)t^{2\delta+1}|^{2}~dt\right)^{1/2}\left(\frac{1}{R_{j}}\int^{R_{j}}_{0}|\mathfrak{B}^{\delta}_{t}(F_1)(x)|^{2}~dt\right)^{1/2}.
\end{eqnarray*}
 Observe that boundedness of  $S^{R,t}_{j,\beta}f_{1}$  is already known. On the other hand, we use induction hypothesis to get the desired boundedness of $\mathfrak{B}^{\delta}_{t}(F_1)$. Indeed, 
 for some $d>1$, we rewrite 
 \begin{eqnarray}
 \mathfrak{B}^{\delta}_{t}(F_1)=\sum_{1\leq k\leq d}(\mathfrak{B}^{\delta+k-1}_{t}(F_1)-\mathfrak{B}^{\delta+k}_{t}(F_1))+\mathfrak{B}^{\delta+d}_{t}(F_1),
 \end{eqnarray}
 where $d$ is chosen such that $\delta+d>n(m-1)/2-1/2$. Now we can write 
 \begin{eqnarray*}
 (\frac{1}{R_{j}}\int^{R_{j}}_{0}|\mathfrak{B}^{\delta}_{t}(F_1)(x)|^{2}~dt)^{1/2}&\leq& \sum_{1\leq k\leq d}\Big(\frac{1}{R_{j}}\int^{R_{j}}_{0}|\mathfrak{B}^{\delta+k}_{t}(F_1)(x)-\mathfrak{B}^{\delta+k-1}_{t}(F_1)(x)|^{2}~dt\Big)^{1/2}\\
 &+& \Big(\frac{1}{R_{j}}\int^{R_{j}}_{0}|\mathfrak{B}^{\delta+d}_{t}(F_1)(x)|^{2}~dt\Big)^{1/2}.
 \end{eqnarray*} 
Observe that for $\delta+d>n(m-1)/2-1/2$, the $(m-1)$-linear Bochner--Riesz means $\mathfrak{B}^{\delta+d}_{t}$ is pointwise dominated by the $(m-1)$-linear Hardy--Littlewood maximal operator $\mathcal{M}_{m-1}$. Hence the operator $F_1\rightarrow \sup_{R>0}(\frac{1}{R_{j}}\int^{R_{j}}_{0}|\mathcal{B}^{\delta+d}_{t}(F_1)(x)|^{2}~dt)^{1/2}$ is bounded from $\prod^{m}_{i=2}L^{p_{i}}$ to $L^{\tilde{p}_{1}}$ for all $1< p_{2},p_{3},\cdots,p_m\leq\infty$ and $\sum^{m}_{i=2}\frac{1}{p_{i}}=\frac{1}{\tilde{p}_{1}}$. The remaining part $\sup_{R>0}(\frac{1}{R_{j}}\int^{R_{j}}_{0}|\mathfrak{B}^{\delta+k}_{t}(F_1)(x)-\mathfrak{B}^{\delta+k-1}_{t}(F_1)(x)|^{2}~dt)^{1/2}$ is dominated by the $(m-1)$-linear Bochner--Riesz square function $(\int^{\infty}_{0}|\mathfrak{B}^{\delta+k}_{t}(F_1)(x)-\mathfrak{B}^{\delta+k-1}_{t}(F_1)(x)|^{2}~\frac{dt}{t})^{1/2}$ with index $\delta+k-1$. Hence using induction hypothesis we get the desired estimates.

 \subsection{Boundedness of $\mathfrak{B}^{\alpha}_{*,0}$}
 We will imitate the approach of Section \ref{B0alpha}. 
 Using the similar decomposition  we get 
 \[m_{\alpha,0}(\vec{\xi}):=m_{\alpha,0,1}(\vec{\xi})+m_{\alpha,0,0}(\vec{\xi}),\]
 where 
 \[m_{\alpha,0,1}(\vec{\xi})=\sum_{j\geq2}\psi_{0}(\frac{|\xi_1|^2}{R^{2}})\psi\Big(2^{j}(1-\frac{|\xi_2|^2}{R^2})\Big)\varphi_{R}(\xi_2)^{\alpha} \Big(1-\frac{|\eta_2|^2\varphi^{-1}_{R}(\xi_2)}{R^2}\Big)^{\alpha}_{+}.\]
 Note that here $\eta_2=(\xi_1,\xi_3,\cdots,\xi_m)$, and we denote the operator corresponding to the above multiplier by $\sum_{j\geq2}\mathfrak{I}^{\alpha,1}_{j,*}$. Hence, invoking Stein and Weiss's identity and applying Cauchy--Schwarz inequality followed by a change of variable  $t\to Rt$ we get 

 \begin{eqnarray*}
&&\Vert \mathfrak{T}^{\alpha,1}_{j,*}(F)\Vert_{L^{p}(\mathbb{R}^{n})}\\
&&\lesssim 2^{-j/4}\Vert\sup_{R>0}(\int^{\sqrt{2^{-j+1}}}_{0}|S^{R,t}_{j,\beta}f_{2}(x)t^{2\delta+1}|^{2}~dt)^{1/2}\Vert_{L^{p_{2}}(\mathbb{R}^{n})}\\
&&\times \Vert \sup_{R>0}(\frac{1}{R_{j}}\int^{R_{j}}_{0}|B^{\psi_0}_{R}\mathfrak{B}^{\delta}_{t}(F_2)(x)|^{2}~dt)^{1/2}\Vert_{L^{\tilde{p}_{2}}(\mathbb{R}^{n})},
\end{eqnarray*} 
where $\widehat{B^{\psi_0}_{R}\mathfrak{B}^{\delta}_{t}(F_2)}(\eta_2)=\psi_{0}(\frac{|\xi_{1}|^{2}}{R^{2}})\widehat{\mathfrak{B}^{\delta}_{t}(F_2)}(\eta_2)$. Note that $\mathfrak{B}^{\delta}_{t}$ is the $(m-1)$ linear Bochner--Riesz means acting on $F_2=(f_1,f_3,\cdots,f_m)$ with index $\delta$.
Now consider $d>1$ such that $\delta+d>n(m-1)/2-{1}/{2}$. Then we write 
\begin{eqnarray*}
 B^{\psi_0}_{R}\mathfrak{B}^{\delta}_{t}(F_2)=\sum_{1\leq k\leq d}(B^{\psi_0}_{R}\mathfrak{B}^{\delta+k-1}_{t}(F_2)-B^{\psi_0}_{R}\mathfrak{B}^{\delta+k}_{t}(F_2))+B^{\psi_0}_{R}\mathfrak{B}^{\delta+d}_{t}(F_2).
\end{eqnarray*}
Now applying Cauchy--Schwarz inequality we get
\begin{eqnarray*}
	&&\Big(\frac{1}{R_{j}}\int^{R_{j}}_{0}|B^{\psi_0}_{R}\mathfrak{B}^{\delta}_{t}(F_2)(x)|^{2}dt\Big)^{1/2}\\
	&&\leq \sum_{1\leq k\leq d}\Big(\frac{1}{R_{j}}\int^{R_{j}}_{0}|B^{\psi_0}_{R}\mathfrak{B}^{\delta+k}_{t}(F_2)(x)-B^{\psi_0}_{R}\mathfrak{B}^{\delta+k-1}_{t}(F_2)(x)|^{2}dt\Big)^{1/2}\\
	&&+ \Big(\frac{1}{R_{j}}\int^{R_{j}}_{0}|B^{\psi_0}_{R}\mathfrak{B}^{\delta+d}_{t}(F_2)(x)|^{2}~dt\Big)^{1/2}.
\end{eqnarray*}
Since $\delta+d>n(m-1)/2-{1}/{2}$, the operator $\mathfrak{B}^{\delta+d}_{t}$ is dominated by the $(m-1)$-linear Hardy--Littlewood maximal function $\mathcal{M}_{m-1}$. Hence we get the desired estimates. On the other hand, boundedness of \[\sum_{1\leq k\leq d}(\frac{1}{R_{j}}\int^{R_{j}}_{0}|B^{\psi_0}_{R}\mathfrak{B}^{\delta+k}_{t}(F_2)(x)-B^{\psi_0}_{R}\mathfrak{B}^{\delta+k-1}_{t}(F_2)(x)|^{2}dt)^{1/2}\] follows in a similar way as in the Subsection \ref{f2f3}. Indeed, this operator will be dominated by the $(m-1)$-linear Bochner--Riesz square function with index $(\delta+k-1)$. Hence from the induction hypothesis we get the desired estimates.
%Note that, similar to the proof of boundedness of trilinear Bochner--Riesz square function, $L^{p}$ boundedness of the $(m-1)$-linear Bochner--Riesz square function follows using the $L^{p}$ boundedness of the $(m-1)$-linear maximal Bochner--Riesz means.
\subsection{Boundedness corresponding to the multiplier $m_{\alpha,0,0}(\vec{\xi})$}
Like the trilinear case, we decompose this multiplier w.r.t the $\xi_2$ variable first and get
\[m_{\alpha,0,0}(\vec{\xi})=m^{1}_{\alpha,0,0}(\vec{\xi})+m^{2}_{\alpha,0,0}(\vec{\xi}).\]
Further, we keep decomposing the multiplier  $m^{1}_{\alpha,0,0}$ w.r.t $\xi_3, \xi_4,\cdots,\xi_m$ variables. Observe that all the intermediate multiplier operators will be handled using the induction hypothesis. Indeed, decomposing w.r.t. the $\xi_m$ variable and applying Stein and Weiss's identity we get that the corresponding operator will be dominated by 
\begin{eqnarray*}
&&\lesssim 2^{-j/4}\Vert\sup_{R>0}(\int^{\sqrt{2^{-j+1}}}_{0}|S^{R,t}_{j,\beta}f_{m}(x)t^{2\delta+1}|^{2}~dt)^{1/2}\Vert_{L^{p_{m}}(\mathbb{R}^{n})}\\
&&\times \Vert \sup_{R>0}(\frac{1}{R_{j}}\int^{R_{j}}_{0}|B^{\psi_{0}}_{R}B^{\psi^{1}_{0}}_{R}\cdots B^{\psi_{0}}_{R}\mathfrak{B}^{\delta}_{t}(F_m)(x)|^{2}~dt)^{1/2}\Vert_{L^{\tilde{p}_{m}}(\mathbb{R}^{n})},
\end{eqnarray*} 
where $\widehat{B^{\psi_{0}}_{R}B^{\psi^{1}_{0}}_{R}\cdots B^{\psi_{0}}_{R}\mathfrak{B}^{\delta}_{t}(F_m)}(\eta_m)=\psi_{0}(\frac{|\xi_{1}|^{2}}{R^{2}})\psi^{1}_{0}(\frac{|\xi_{2}|^{2}}{R^{2}})\cdots \psi_{0}(\frac{|\xi_{m}|^{2}}{R^{2}})\widehat{\mathfrak{B}^{\delta}_{t}(F_m)}(\eta_m)$. Note that $\mathfrak{B}^{\delta}_{t}$ is the $(m-1)$ linear Bochner--Riesz means acting on $F_m=(f_1,f_2,\cdots,f_{m-1})$ with index $\delta$.  Boundedness of this operator can be proved using similar ideas as the previous subsection. Similarly, boundedness of the operator corresponding to the multiplier $m^{2}_{\alpha,0,0}$ can be completed using decompositions w.r.t. the variables $\xi_3,\xi_4,\cdots,\xi_m$  along with the Stein and Weiss's identity and induction hypothesis. We skip this part to avoid repetitions of the previous methods. Boundedness of $\mathfrak{B}^{\alpha}_{*}$ in dimension $n=1$ can be completed in a similar way.
\end{proof}
\subsection{Proof of Remark \ref{rem}}
Boundedness of $\mathcal{G}^{\alpha}$ can be proved in a similar manner as the trilinear Bochner--Riesz square function  using Theorem \ref{mlinear} and \ref{mlinear1}.  
 %\begin{Theorem}
 	
 %\end{Theorem}

\vspace{.25in}
\noindent
{\bf Acknowledgements:} 
The author extends sincere gratitude to Prof. Saurabh Shrivastava for his constant encouragement to pursue independent research and for providing valuable feedback on earlier versions of this draft. The author also wishes to express deep appreciation to previous collaborators on this research topic, Dr. Jotsaroop Kaur and Dr. Surjeet Singh Choudhary. This work is supported by DST Inspire Faculty Award (Registration no. IFA23-MA 191).

\vspace{.3in}
%\newpage

\end{document}